\documentclass{article}
\usepackage{graphicx,color}

\usepackage{amsmath,amssymb,amsfonts,epsfig, %authblk,
mathrsfs,hyperref,doi}
\usepackage[shortlabels]{enumitem}
\usepackage{ulem} %new
\numberwithin{equation}{section}

\newtheorem{prop}{Proposition}
\numberwithin{prop}{section}

\newtheorem{defi}[prop]{Definition}
\newtheorem{cor}[prop]{Corollary}

\newtheorem{thm}[prop]{Theorem}

\newcommand{\ex}{\medskip \noindent {\bf Example }}

\newcommand{\proof}[1]{\vspace{1.5ex}\noindent{{\bf Proof:} #1 
}
\vspace{1.5ex}
}

\newcommand{\beqa}{\begin{eqnarray*}}
\newcommand{\eeqa}{\end{eqnarray*}}

\newcommand{\beqan}{\begin{eqnarray}}
\newcommand{\eeqan}{\end{eqnarray}}

\newcommand{\beq}{\begin{equation}}
\newcommand{\eeq}{\end{equation}}

\newcommand{\tsf}[2]{\protect{{\textstyle{ \frac{#1}{#2}}}}}

\newcommand{\mc}[1]{\protect{{\mathscr{#1}}}}

\newcommand{\WFG}{\protect{\mathrm{WF_G}}}

\newcommand{\Z}{\mathbb Z}
\newcommand{\R}{\mathbb R}
\newcommand{\Sp}{\mathbb S}
\newcommand{\B}{\mathbb B}

\title{The $G$-wavefront set and the twisted convolution product}

\author{Dorothea Bahns\footnote{Mathematisches Institut, Georg-August-Universit\"at G\"ottingen, Bunsenstr. 3-5, D - 37073 G\"ottingen, Germany,  dbahns@mathematik.uni-goettingen.de} \  and Ren\'e Schulz\footnote{Institut f\"ur Analysis, Leibniz Universit\"at Hannover, Welfengarten 1, D - 
30167 Hannover, Germany, rschulz@math.uni-hannover.de}}

\date{ }

\begin{document}

\maketitle

\begin{abstract}
\noindent
We give a sufficient criterion for the existence of the twisted convolution product of two tempered distributions as a tempered distribution, and we list examples of algebras with respect to this and related products  contained in $\mc S^\prime$.

\end{abstract}

%%%%%%%%%%%%%%%%

\section{Introduction}\label{sec:intro}
\setcounter{equation}{0}

Phase space quantum mechanics has a long history. Dating back to Wigner, Groenewald, Moyal and others, it comprises the attempt to study quantum mechanics not in terms of (unbounded) 
operators on a Hilbert space, e.g. position and momentum operator on $L^2(\R^m)$, but by considering instead classical functions on phase space 
$T^*\R^m$ endowed with a noncommutative product, see e.g.~\cite{wigner,groenewald, moyal}. For two Schwartz functions $f, g \in \mc S(\R^{2m})$, such a product is the twisted convolution 
\beq\label{eq:twConv}
f\star g (x)=\int_{\R^{2m}}  f(x-y)g(y)e^{-\frac i  2 x^T (\theta y)}\,dy \ ,
\eeq
where $\theta$ is the canonical symplectic form on $\R^{2m}$, $\theta=\left(\begin{array}{cc}0&1_m\\-1_m&0\end{array}\right)$ and where $x=(q,p)\in T^*\R^m$ is interpreted as a phase space variable. Depending on context, the following, similar product (``twisted convolution {\sl product}'') 
\beq\label{eq:twConvP}
f\ast g (x)=\mc F^{-1}_{\xi \rightarrow x} \int_{\R^{2m}} \mc F (f)(\xi-\eta) \, \mc F (g)(\eta)\, e^{-\frac i 2 \xi^T (\theta \eta)} \,d \eta \ ,
\eeq
where $\mc F$ denotes the Fourier transform, also appears in the literature. 
Both formulae obviously again yield  a Schwartz function for any matrix $\theta$ and $f,g \in \mc S$, and the products are associative for any  antisymmetric matrix.
From a purely mathematical viewpoint, the product \eqref{eq:twConvP} arises as the symbol product in the Weyl calculus of pseudo-differential operators and of course, it has been used in different contexts from the early days of quantum mechanics, see e.g~\cite{neumann}.
For $\theta=0$, it corresponds to the pointwise product.
Some 20 years ago, phase space quantum mechanics experienced a revival in theoretical physics in the context of noncommutative geometry. In fact, one of the basic examples of a noncommutative spacetime -- dating back to an even earlier time -- is given by postulating canonical commutation relations for the coordinate operators~\cite{DFR} and makes use of the above products in suitable spaces of functions instead of studying these operators.\footnote{In the study of quantum field theory on such a noncommutative space, it was even customary to simply put quadratic phases into all products in the (Euclidean) Feynman rules until it was dicovered that these rules had to be modified in a more substantial way in order to have a well-defined theory~\cite{DFR, BDFP_unitary}.}. In mathematics, this approach was viewed  as an early  example of  a Rieffel deformation~\cite{rieffel_onDFR}, a method to deform a  commutative algebra into a noncommutative one.

In studying quantum field theory on such a simple noncommutative spacetime, the question arose how to extend the above products (viz. the Weyl calculus) to (tempered) distributions. This question, which is  interesting in its own right, had been raised and partially answered in different contexts, e.g. by Maillard~\cite{maillard} as well as Gracia-Bond\'ia and Varilly~\cite{gbVar_moyal1, gbVar_moyal2} (in the context of quantum mechanics), by Dubois-Violette, Kriegl, Maeda and Michor~\cite{michor} as well as Lechner and Waldmann~\cite{lw} (in the context of noncommutative geometry) and by others, among them Zahn~\cite{zahn} (in the context of quantum field theory on noncommutative spaces, see also~\cite{B_UVIR}). 
More generally, the question can be raised in the context of the Weyl calculus (see e.g. Folland's~\cite{folland}): Every continuous operator from $\mc S(\R^n)\rightarrow \mc S'(\R^n)$ has a Weyl symbol in $\mc S'(\R^{2n})$, and hence the composability of two operators is equivalent to the existence of the Weyl product of their symbols.

Using the global or Gabor wavefront set, which is a microlocal concept adapted to tempered distributions and in particular measures the behaviour at infinity (ignored by the classical wavefront set), we give a sufficient criterion that guarantees the existence of the twisted convolution product that was first published in the doctoral thesis~\cite{schulz_thesis}. To give this result a broader audience and to embed it into context, we present this result and its proof in more detail here, give examples of topological algebras of tempered distributions with respect to the twisted convolution and the twisted convolution product and compare our result to the literature.

\section{The $G$-wavefront set}

Originally introduced by H\"ormander~\cite{hoerm_WFG} to study the propagation of of singularities of quadratic hyperbolic operators, the $G$-wavefront set $\WFG$ ($G$ for global or Gabor) is a microlocal tool that takes into account the behaviour of a (tempered) distribution at infinity. While the classical wavefont set measures the behaviour of a distribution in cones in the cotangent space $T^*\R^n$, the $G$-wavefront set is characterized by open subsets 
in the boundary $\Sp^{2n-1} \simeq (\R^{2n}\setminus \{0\}/\sim)$, with $\sim$ given by the dilation, of the compactification\footnote{Taking the compactification $\B^{n}\times \B^n$ of $T^*\R^n\simeq \R^n \times \R^n$,  one derives the so-called $SG$-wavefront set. A detailed comparison and analysis of these  wavefront sets can be found in~\cite{schulz_thesis} where it is also proved that the $G$-wavefront set coincides with Nakamura's homogeneous wavefront set.} $\B^{2n}$ of the cotangent space $
 T^* \R^n\simeq \R^{2n}$. For the reader's convenience, we recall here the basic definition and properties. We choose the characterization of $\WFG$  in terms of short term Fourier transform, instead of the microlocal one, as it requires the least  prior knowledge. Its equivalence to H\"ormander's definition was proved in~\cite{rodWahlb}. 
Let $\psi \in \mc S(\R^n)\setminus \{0\}$ be a window  function, e.g. the (normalized) Gaussian. Then the short term Fourier transform of $f\in \mc S(\R^n)$ is 
\beq
\mc V_\psi(f)(x,\xi)= \tsf 1 {(2\pi)^{n/2}\, \|\psi\|_2}\int_{\R^n} f(y)\overline{\psi(y-x)} \,e^{-i\xi \cdot x} \, dy
\eeq  
This gives a continuous map $\mc V_\psi : \mc S (\R^n) \rightarrow \mc S(\R^{2n})$ which extends to a map $\mc V_\psi : \mc S^\prime (\R^n) \rightarrow \mc S^\prime(\R^{2n})$ by
\[
\mc V_\psi(u)(x,\xi)=\tsf 1 {(2\pi)^{n/2}\, \|\psi\|_2} (u,\mc M_\xi \tau_x \psi)
\]
where $\mc M_\xi g(x) = e^{i\xi \cdot x}g(x)$ and $\tau_x$ denotes the translation and where the pairing is $(u,g)=u(\overline g)$. Observe that for $u \in \mc S^\prime$, $(x,\xi)\mapsto \mc V_\psi(u)(x,\xi)$ is a polynomially bounded continuous function.
\begin{defi}\label{def:WFG}
Let $u \in \mc S^\prime(\R^n)$, let $(x_0,\xi_0) \in \R^{2n}\setminus \{0\}$. Then $(x_0,\xi_0)$ is \emph{not} in the  wavefront set $\WFG(u)$ of $u$ if and only if for one (hence for all) $\psi \in \mc S \setminus \{0\}$, there is an open cone $\Gamma_0 \subset \R^{2n}\setminus \{0\}$ containing $(x_0,\xi_0)$ such that
\[
(1+\|(x,\xi)\|^2)^k \, |\mc V_\psi u(x,\xi)| 
\]    
is bounded on $\Gamma_0$ for all $k=0,1,2,\dots$.
\end{defi}

Observe that $\WFG$ is a closed conic (w.r.t. variable and covariable jointly) subset of $\R^{2n}\setminus \{0\}$. Since the cones in the above definition  are localized ``at infinity" (i.e. in the compactification's boundary), all local information is lost, e.g. the $\delta$-distribution at $a \in \R^n$ always has wavefront set $\WFG(\delta_a)=\{0\}\times (\R^n\setminus \{0\})$ for any $a$. In fact, lack of smoothness at finite points leads to a contribution $\{0\}\times (\R^n\setminus \{0\})$ to the wavefront set, lack of decay and linear oscillations at infinity are encoded in $(\R^n\setminus \{0\})\times \{0\}$, while the rest of $\R^{2n}\setminus \{0\}$ encodes higher than linear oscillations at infinity (see~\cite{rodWahlb, schulz_thesis} and the discussion below).

\begin{prop}[Properties of $\WFG$ {[H\"ormander]}] \label{prop:properties}Let $u \in \mc S^\prime(\R^n)$. We have
\begin{enumerate}[(i)]
\item \label{reg} Regularity: $u \in  \mc S \subset \mc S^\prime$ $\Leftrightarrow$ $\WFG(u) = \emptyset$.
\item \label{Four} Symmetry under Fourier transform: $(x,\xi) \in \WFG(u)$ $\Leftrightarrow$ $(\xi,-x) \in \WFG(\mc F  u)$ 
\item \label{chirp} Multiplication by a chirp: For $A$ a real, symmetric $n \times n$ matrix, we have
\[
(x,\xi) \in \WFG(u) \quad \Leftrightarrow \quad (x, \xi +Ax) \in \WFG(e^{\frac i  2 x^T (Ax)}u)
\] 
\item \label{pullb} Pullback: Let $A: \R^n\rightarrow \R^m$ be linear. If $\WFG \cap \{(0,\xi)| A^T\xi =0\}= \emptyset$, the pullback $A^*u$ of $u$ along $A$ is uniquely defined as a tempered distribution and 
\[
\WFG(A^*u) \subset \{(x,A^T \xi) \,| \, (Ax,\xi)  \in \WFG(u)\}\, \cup \, \ker A \times \{0\}
\]
\end{enumerate}
\end{prop}

\textbf{Proof:} Statement \ref{reg} is proved in \cite{hoerm_WFG}, statements \ref{Four} and \ref{chirp} in \cite[Thm 18.5.9]{hoerm_3}. Statement  \ref{pullb} is proved in \cite{hoerm_WFG} (Proposition 2.3 and 2.9). Observe, however, that in the second of these propositions,  the contribution from the kernel of $A$ was accidentally omitted.

As an easy corollary observe that from \ref{Four} above and the wavefront set of $\delta_a$ it follows that $\WFG(e^{ia \cdot})=(\R^n\setminus \{0\})\times \{0\}$ for any $a \in \R^n$.

The $G$-wavefront set is appropriate to test whether the pointwise product of two tenpered distributions exists and is again temepered:
\begin{thm}\label{thm:pointwise}
Let $u, v \in \mc S (\R^n)$ be such that $(0,\xi) \in \WFG(u)$ implies that $(0,-\xi) \not \in \WFG (v)$. Then the pointwise product of $u$ and $v$ is uniquely defined as the extension of the pointwise product of two Schwartz functions and it is again tempered.  It is given by the pullback $\mathrm{diag}^*(u\otimes v)$ of the tensor product  along the diagonal map $\mathrm{diag}(x)=(x,x)$.
\end{thm}

The proof was not explicitly given in~\cite{hoerm_WFG} which concentrated instead on pairings of distributions. It follows immediately from our more general theorem below. 

%%%%%%%%%%%%%%

\section{The twisted convolution product on $\mc S^\prime$}

To extend the twisted convolution product to tempered distributions, we rewrite the product of  $f, g \in \mc S(\R^{2m})$ in terms of the tensor product and the pullback along the diagonal map,
\beqa
f\ast g (x) &=& \left(\mc F^{-1}_{(k,p) \rightarrow (x,y)} \left( \mc F (f)(k) \otimes \mc F (g)(p)\, e^{-\frac i 2 k^T (\theta p)} \,\right)\  \right)_{x=y}
\\
&=& \mathrm{diag}^*\left(\mc F^{-1}  \left( \mc F (f)(k) \otimes \mc F (g)(p)\, e^{-\frac i 2 k^T (\theta p)} \,\right)\  \right)(x)
\\
&=& \mathrm{diag}^*\left(\mc F^{-1}  \left( (\mc F (f) \otimes \mc F (g))\, \cdot \, \chi_\theta \,\right) \right)(x)
\eeqa
where the equality of the first line with \eqref{eq:twConvP} follows with a simple variable transform ($k \rightarrow k-p$) and where $\chi_{\theta}$ is the chirp on $\R^{2m} \times  \R^{2m}$,
\[
\chi_\theta (K) = \exp(-\tsf i 2 K^T (\Theta K) )\ , \qquad \Theta=\tsf 1 2 \left(\begin{array}{cc} 0 & \theta \\ -\theta & 0\end{array}\right)
\]
since for $K=\left(k \atop p\right)$, we have $K^T (\Theta K)=\frac 1 2 (k^T\theta p - p^T \theta k)= k^T\theta p$. Observe that indeed, we have $\Theta^T=\Theta$. 

We now prove our main result.
\begin{thm} 
Let $\theta$ be an antisymmetric $n\times n$ matrix, and let $u, v \in \mc S^\prime (\R^{n})$ be such that 
\beq\label{eq:cond}
x - \tsf 1 2 \, \theta \xi = 0 
\eeq 
has \emph{no} solution $(x, \xi)$ with $(x,\xi) \in \WFG(u)$ and $(x,-\xi) \in \WFG(v)$.  Then the twisted convolution product of $u$ and $v$ is defined as
\[
u\ast_\theta v = \mathrm{diag}^*\left(\mc F^{-1}  \left( (\mc F (u)  \otimes \mc F (v)) \, \chi_\theta  \,\right)\  \right)
\]
and its $G$-wavefront set $\WFG(u\ast_\theta v)$ is contained in 
\[
\Big \{(x+\tsf  1 2 \theta\eta,\xi+\eta)\, \Big| \,  (x,\xi) \in \WFG(u) \cup \{0\}\, , \, (y,\eta)\in \WFG(v)\cup \{0\} \, ,  y %-\tsf 1 2 \theta \xi
=x+\tsf 1 2 \theta( \eta + \xi)
\, 
\Big\}\setminus \{0\}
\]
It yields a sequentially continuous map $\mc S^\prime_{\WFG(u)} \times  \mc S^\prime_{\WFG(v)} \rightarrow \mc S^\prime_{\Gamma}$.  
\end{thm}

Here, for any closed cone $\Gamma \subset \R^{2n}\setminus \{0\}$, $\mc S^\prime_{\Gamma}(\R^n)$ is the space of all tempered distributions with $\WFG$ contained in $\Gamma$, endowed with the H\"ormander topology.

\proof{From Proposition~\ref{prop:properties} \ref{Four}, we deduce that the wavefront set of $\mc Fu \otimes \mc F v$ is contained in 
\[
\big\{(\xi,\eta; -x,-y) \ | \  (x,\xi) \in \WFG(u) \cup \{0\} \mbox{ and } (y,\eta)\in \WFG(v)\cup \{0\}\, 
\big\}\setminus \{0\}
\] 
and from Proposition~\ref{prop:properties} \ref{chirp}, it follows that the wavefront set of $(\mc Fu \otimes \mc F v) \ \cdot \chi_\theta$ is contained in 
\beqa
&&\big\{\Big(\left( \xi \atop \eta\right) ; \underbrace{-\left(x \atop y\right) - \Theta \left( \xi \atop \eta\right)}_{\displaystyle =\left({-x - \tsf 1 2 \theta \eta} \atop { - y + \tsf 1 2 \theta \xi}\right)} \Big) \ | \  (x,\xi) \in \WFG(u) \cup \{0\} \mbox{ and } (y,\eta)\in \WFG(v)\cup \{0\}\, 
\big\}\setminus \{0\}
\eeqa
Taking the inverse Fourier transform then yields (again by Proposition~\ref{prop:properties} \ref{Four})
\[
\left\{\left(x + \tsf 1 2 \theta \eta,  y - \tsf 1 2 \theta \xi; \xi, \eta\right) \ | \  (x,\xi) \in \WFG(u) \cup \{0\} \mbox{ and } (y,\eta)\in \WFG(v)\cup \{0\}\, 
\right\}\setminus \{0\}
\]
for the $G$-wavefront set of $\mc F^{-1}((\mc Fu \otimes \mc F v) \cdot \chi_\theta)$. By Proposition~\ref{prop:properties} \ref{pullb},  the pullback onto the diagonal is well-defined if the above has empty intersection with 
\[
\{ (0;\xi, \eta )\ | \ \xi +\eta =0 \} \ .
\]
This is guaranteed by our assumption on the wavefront sets of $u$ and $v$: If $x+ \tsf 1 2 \theta \eta = 0 =  y - \tsf 1 2 \theta \xi$  and $\eta=-\xi$, then 
\[
x- \tsf 1 2 \theta \xi = 0\ , \quad (y,\eta)=(x,-\xi)
\]
and by assumption, the first equation does not have a solution $(x,\xi) \in \WFG(u)$ such that $(x,-\xi) \in \WFG(v)$.

By Proposition~\ref{prop:properties} \ref{pullb}, its wavefront set is contained in 
$$\{(x,\xi+\eta)\ | \ (x,x;\xi,\eta) \in \WFG(\mc F^{-1}((\mc Fu \otimes \mc F v) \cdot \chi_\theta)) \}\  \cup \ \ker (\mathrm{diag}) \times \{0\}
$$ 
which in turn is equal to
\beqa
\Big \{(x+\tsf  1 2 \theta\eta,\xi+\eta)\, \Big| \,  (x,\xi) \in \WFG(u) \cup \{0\}\, , \, (y,\eta)\in \WFG(v)\cup \{0\} \, ,  y 
=x+\tsf 1 2 \theta( \eta + \xi)
\, 
\Big\}\setminus \{0\}
\eeqa
where we got rid of the total zero $\ker(\mathrm{diag}) \times \{0\}=\{0\} \in \R^{2n}$, as it is by definition not contained in the wavefront set.

The continuity statement follows from the sequential continuity of the elementary operations above which constitute the twisted convolution product.  \hfill$\square$

}

The statement in Theorem~\ref{thm:pointwise} for the pointwise product now follows automatically, by setting $\theta=0$. For the other extreme, where  $\theta$ is invertible (hence $n$ must be even), the condition on the wavefront sets in the theorem reads
\[
(x, 2 \, \theta^{-1} x) \in \WFG(u) \quad \mbox{ implies }\quad  (x, -2 \, \theta^{-1} x) \not \in \WFG(v) \ .
\]

Observe that H\"ormander's classical wavefront set would be unsuitable for this investigation, as -- but for homogeneous distributions --  it lacks the symmetric behaviour under Fourier transform which $\WFG$ possesses.

%%%%%%%%%%%%%%

\section{Algebras in $\mc S^\prime$}

An immediate consequence of our main theorem is the following:

\begin{cor} 
Let $\theta$ be an antisymmetric $n\times n$ matrix. Let  $\Gamma$ be a conic subset of $\R^{2n} \setminus 0$ such that $(x,\xi) \in \Gamma$ implies that $(x,-\xi) \not\in \Gamma$, then pointwise multiplication and twisted convolution product of any two $u,v \in \mc S^\prime_\Gamma(\R^n)$ yield well-defined tempered distributions.
\end{cor}

Observe however, that the wavefront set of the resulting distributions will in general no longer be an element of $\mc S_\Gamma^\prime$.

\begin{cor}\label{cor:schulz} 
Let   $\theta$ be an antisymmetric $n\times n$ matrix. Let  $\Gamma_2$ be a conic subset of $\R^n \setminus 0$ that is closed under addition and let $\Gamma_1$ be a conic subset of $\R^n$ such that 
\beq\label{eq:corSchulz}
x \in \Gamma_1, \xi \in \Gamma_2 \quad \mbox{ implies } \quad x +\tsf 1 2 \theta \xi \in \Gamma_1 
\eeq
Then $\mc S_{\Gamma_1\times \Gamma_2}^\prime(\R^n)$ is an algebra under pointwise multiplication and under the twisted convolution product.

\end{cor}

\proof{
Since $\Gamma_2$ is closed under addition, but does not contain 0, condition \eqref{eq:cond} in the theorem is trivially satisfied. The pointwise product has wavefront set contained in $\Gamma_1\times \Gamma_2$, since $\Gamma_2$ is closed under addition. This is also true for the twisted convolution product (with $\theta$ non trivial) since by definition, $x + \frac 1 2 \theta \eta \in \Gamma_1$ for any $\eta \in \Gamma_2$ and $x \in \Gamma_1$.\hfill $\square$

}

\ex 
Let $n=2$  and let $\theta=\left(\begin{array}{cc}0&-1\\1&0\end{array}\right)$. Let  $\Gamma_2$ be  the positive lightcone $\{(\xi_0,\xi_1)| \xi_0^2\leq \xi_1^2\,, \xi_0\geq 0\}$. Then $\Gamma_2$ is closed under addition, and $\xi \in \Gamma_2$ implies $-\xi \not\in \Gamma_2$. Setting $\Gamma_1$ equal to the left half space will satisfy the condition \eqref{eq:corSchulz}.
\medskip

Obserse that these algebras are non-empty since for any closed closed conic set $\Gamma$ in $\R^{2n}\setminus \{0\}$, there is a distribution $u \in \mc S^\prime_\Gamma(\R^n)$, see e.g.~\cite{schulzWahlb}.

We now consider some older results from  the literature in view of our framework. To this end, we first observe that for $f, g \in \mc S$, the twisted convolution \ref{eq:twConv} can be written as
\beq\label{eq:twCtwCP}
f \star g (y) = \mc F_{x\rightarrow y}( (\mc F^{-1}f \ast \mc F^{-1}g)(x))
\eeq
in terms of the twisted convolution product. We deduce:

\begin{cor}\label{cor:michor} The space $\mc S_{\{0\}\times \R^n\setminus \{0\}}^\prime(\R^n)$ forms an algebra under the twisted convolution~$\star$.
\end{cor}

To see this, let  $u,v \in \mc S_{\{0\}\times \R^n\setminus \{0\}}^\prime(\R^n)$. We conclude from our theorem and the behaviour of $\WFG$ under Fourier transform (Proposition~\ref{prop:properties} \ref{Four}) that the twisted convolution product of $\mc F^{-1}u$ and $\mc F^{-1}v$ is well-defined, since \eqref{eq:cond} has no solution at all with $(x,\xi) \in \R^n\setminus \{0\}\times \{0\}$. We conclude moreover,  that the wavefront set of $\mc F^{-1}u \ast \mc F^{-1}v$ is contained in $\R^n\setminus \{0\} \times \{0\}$, and hence, the wavefront set  we obtain by taking the additional Fourier transform  (as prescribed by \eqref{eq:twCtwCP}) is again $ \{0\}\times \R^n\setminus \{0\} $. \medskip

Consider in this context the following  result from~\cite{michor}: Let $\mathcal O_M^\prime(\R^{n})$ denote the space of speedily decreasing distributions, i.e. the dual space of $\mathcal O_M(\R^{n})$, 
\[
\{ f \in C^\infty(\R^n) | \mbox{ for any } \alpha \in \mathbb N_0^n \mbox{ there is } k \in \Z \mbox{ s.t. } (1+|x|^2)^k \partial^\alpha f \mbox{ is bbd} 
\}\, ,
\] 
endowed with the topology given by all these seminorms. 
Then for $u,v \in \mathcal O_M^\prime(\R^{2m})$, the twisted convolution \eqref{eq:twConv} $u\star v$ with the standard symplectic $\theta$ is again in  $\mathcal O_M^\prime(\R^{2m})$.

To see the relation with Corollary~\ref{cor:michor}, we first observe  with~\cite{michor} that
\[
\mc S \subset \mathcal O_C \subset \mathcal O_M \subset \mc S^\prime
\quad \mbox{ and }\quad
\mc S \subset \mathcal O_M^\prime \subset \mathcal O_C^\prime \subset \mc S^\prime
\]
where $\mathcal O_C(\R^n)$ is the space of smooth functions such that there is $k \in \Z$ s.t. $(1+|x|^2)^k \partial^\alpha f$ is bounded for all  multindices $\alpha$, and that moreover,  by Fourier transform,  
$\mathcal O_M^\prime(\R^{n}) \simeq \mathcal O_C(\R^n)$. 

It is easy to see that the $G$-wavefront set (in the sense of definition~\ref{def:WFG}) of $a \in \mathcal O_C(\R^n)$ is contained in $\R^n\setminus \{0\} \times \{0\}$, as the polynomial growth in $\|(x,\xi)\|$ can (for $x\neq 0$) be estimated by a linear combination of derivatives (w.r.t. $x$) of $\psi(y-x)\,e^{i x\xi}$, and  rewriting
\[
 \partial_x^\beta \left(\psi(y-x)\,e^{-i\xi x}\right) \, = \,  
e^{-i\xi y}\, \partial_x^\beta \left(\psi(y-x)\,e^{-i\xi(x-y)}\right) \ , 
\]
the claim follows by integration by parts, using that  $y\mapsto a(y)\,e^{-i\xi y}$ is in $\mathcal O_C(\R^n)$ for $a \in \mathcal O_C(\R^n)$. From the isomorphism $\mathcal O_M^\prime(\R^{n}) \simeq \mathcal O_C(\R^n)$ and the behaviour of $\WFG$ under Fourier transform (Proposition~\ref{prop:properties}\ref{Four}) we conclude that the wavefront set of $u \in \mathcal O_M^\prime$ is contained in $\{0\} \times \R^n\setminus \{0\}$.

Hence, we now obtain from corollary~\ref{cor:michor}, that the twisted convolution of $u,v \in \mathcal O_M^\prime$ is well defined, and its $G$-wavefront set is again contained in $\{0\} \times \R^n\setminus \{0\}$. Observe that if $\mathcal O_C$ were equal to $\{u \in \mc S^\prime(\R^n) \, | \, \WFG(u) \subset \R^n\setminus \{0\} \times \{0\} \}$, then we would get exactly the result of~\cite{michor} that the twisted convolution $u\star v$ is again in $\mathcal O_M^\prime$.

On the other hand, we cannot mimick the construction from~\cite{gbVar_moyal1, gbVar_moyal2} in our context. There,  the so-called Moyal algebra was given as the intersection of 
\[
\mc M_L=\{ S \in \mc S^\prime | S \times f \in \mc S \mbox{ for all } f \in \mc S \}\quad \mbox{ and } \quad
\mc M_R=\{ S \in \mc S^\prime | f \times S \in \mc S \mbox{ for all } f \in \mc S \}
\] 
with yet another (related) twisted product $\times$ given in terms of the standard symplectic form $\theta$. If, in our context, we were to replace this product  by  the twisted convolution product and to consider  
$\mc M_{L,\ast}:=\{ u \in \mc S^\prime | u \ast f \in \mc S \mbox{ for all } f \in \mc S \}$, then we would find that any $u \in \mc S^\prime$ has a well-defined twisted convolution product with any $f \in \mc S$. However,  the inclusion we have for the $G$-wavefront set of $u\ast f$ is too crude to say when exactly the twisted convolution product is a Schwartz function: Asking 
\[
\Big \{(x+\tsf  1 2 \theta\eta,\xi+\eta)\, \Big| \,  (x,\xi) \in \WFG(u) \cup \{0\}\, , \, (y,\eta)\in \WFG(f)\cup \{0\} \, ,  y=x+\tsf 1 2 \theta( \eta + \xi)\Big \} \setminus \{0\}
 \]
for $f \in \mc S$, i.e. $\WFG(f)=\emptyset$, to be empty is too strong, as it amounts to requiring $u$ to be Schwartz itself.\bigskip

Nevertheless, the global wavefront set is a tool that provides a sufficent criterion for the existence and makes it possible to set up spaces of  tempered distributions that form algebras under pointwise multiplication, $\star$ or $\ast$.
It remains an interesting open problem to investigate other (exact) Rieffel deformations in this context, e.g. the ones studied in~\cite{lw,bg}.

%%%%%%%%%%%%%%%%%%%%%%%%%%%%%%%%%%%

\end{document}